\documentclass
[
11pt]
{amsart}

\usepackage{xcolor}
\usepackage[]{hyperref}
\hypersetup{
	colorlinks=true,
	linkcolor=blue,    
	citecolor=green,   
	urlcolor=red       
}
\usepackage{cite}
\usepackage{amssymb}
\usepackage{amsmath}
\usepackage{amsthm}
\usepackage{amsfonts}
\usepackage{bbm}
\usepackage{enumerate} 
\usepackage{graphicx}
\usepackage{xcolor}
\usepackage[toc,page]{appendix}
\usepackage{subfigure}
\usepackage{mathtools}
\usepackage[normalem]{ulem}
\usepackage{comment}
\usepackage{hyperref}

\newcommand{\R}{\mathbb{R}}
\renewcommand{\P}{\mathbb{P}}
\newcommand{\supp}{\mathrm{supp}}
\newcommand{\Per}{\mathrm{Per}}
\newcommand{\E}{\mathbb{E}}
\renewcommand{\ge}{\geqslant}
\renewcommand{\le}{\leqslant}

\newcommand{\BV}{\mathrm{BV}}
\renewcommand{\epsilon}{\varepsilon}
\newcommand{\Exp}{\mathrm{Exp}}

\makeatother
\numberwithin{equation}{section}
\newtheorem{theorem}{Theorem}[section]

\newtheorem{proposition}[theorem]{Proposition}

\newtheorem{lemma}[theorem]{Lemma}

\theoremstyle{definition}
\newtheorem{definition}[theorem]{Definition}

\theoremstyle{remark}
\newtheorem{remark}[theorem]{Remark}
\newtheorem{notation}{Notation}[section]

\makeatletter
\@namedef{subjclassname@2020}{%
	\textup{2020} Mathematics Subject Classification}
\makeatother

\title[Fractional heat content]{Fractional heat content asymptotics for Carnot groups}
\author[Rohan Sarkar]{Rohan Sarkar{$^{\dag }$}}
\address{ Department of Mathematics and Statistics\\
	Binghamton University\\
	Binghamton, NY 13902,  U.S.A.}
\email{rohansarkar20@gmail.com}


\begin{document}
	\begin{abstract}
	We propose a novel approach for studying small-time asymptotics of the fractional heat content of $C^2$ non-characteristic domains in Carnot groups. Denoting the sub-Laplacian operator by $\mathcal{L}$, the fractional heat content of a bounded domain $\Omega$ is defined as $Q^{(\alpha)}_\Omega(t)=\int_{\Omega}u_\alpha(x,t) dx$, where $u_\alpha$ is the solution to the heat equation corresponding to the fractional sub-Laplacian $\mathcal{L}_\alpha:=\mathcal{L}^{\alpha/2}$ with Dirichlet boundary condition on $\Omega$. We prove that for $1\le \alpha\le 2$, there exists explicit rate function $\mu_\alpha: (0,\infty)\to (0,\infty)$ such that
	\begin{align*}
		\lim_{t\to 0}\frac{|\Omega|-Q^{(\alpha)}_\Omega(t)}{\mu_\alpha(t)}=|\partial \Omega|_H,
	\end{align*}
	where $|\Omega|$, $|\partial \Omega|_H$ are the volume and horizontal perimeter of $\Omega$ respectively. Moreover, the rate function $\mu_\alpha$ coincides with the same for the Euclidean case.
	\end{abstract}
	
	\maketitle
	\section{Introduction and the main result}
	We start by defining the fractional heat content of a domain in Euclidean space. Let $\Omega\subset \R^n$ be a bounded domain and $\alpha\in (0,2]$. Consider the following parabolic equation with Dirichlet boundary condition: 
	\begin{align*}
		\left(\frac{\partial}{\partial t}+\Delta_\alpha\right)u_\alpha(t,x)&=0 \quad \mbox{in $(0,\infty)\times\Omega$} \\
		u_\alpha(t,x)&=0 \quad \mbox{in $(0,\infty)\times\partial\Omega$} \\
		u_\alpha(0,x)&=1 \quad\mbox{$\forall x\in\Omega$},
	\end{align*}
	where $\Delta_\alpha$ is the fractional Laplacian defined by
	\begin{align*}
		\Delta_\alpha u(x)=C(n,\alpha) \ \mathrm{p.v.}\int_{\R^n}\frac{u(x)-u(y)}{|x-y|^{n+\alpha}} dy, \quad C(n,\alpha)=\frac{2^\alpha\Gamma(\frac{n+\alpha}{2})}{\pi^{n/2}\Gamma(\frac{\alpha}{2})}.
	\end{align*}
	Under standard regularity conditions on the boundary of $\Omega$ it is known that the solution to above non-local equation is unique.
	The fractional heat content of $\Omega$ is defined as $Q^{(\alpha)}_\Omega(t)=\int_{\Omega}u_\alpha(t,x)dx$. In particular when $\alpha=2$, we recover the classical heat content, and its small-time asymptotics have been of special interest as it depends on the geometry of the domain. The first order asymptotic  in the diffusion case ($\alpha=2$) was originally proved by van den Berg and Davies \cite{BergDavies1989}, and subsequently, van den Berg and LeGall \cite{BergLeGall1994} obtained the following second order expansion for bounded domains with $C^3$ boundary:
	\begin{align}\label{eq:BergDavies}
		Q^{(2)}_\Omega(t)=|\Omega|-\sqrt{\frac{2t}{\pi}}\sigma(\partial\Omega)+\frac{t}{4}\int_{\partial\Omega} H_\Omega d\sigma+\mathrm{o}(t),
	\end{align}
	where $|\Omega|$ is the volume of $\Omega$, $\sigma$ is the surface measure on $\partial\Omega$, and $H_\Omega$ denotes the mean curvature of the boundary points. Later, the above asymptotic formula was extended to Riemannian manifolds by van den Berg and Gilkey \cite{BergGilkey1994}, and an asymptotic expansion holds similar to \eqref{eq:BergDavies} adapted to Riemannian manifolds. We also refer to excellent surveys by Gilkey \cite{Gilkey1999, Gilkey2009} on heat content asymptotics for Riemannian manifolds.
	
	While the results for diffusion case are known for many years, asymptotics of $Q^{(\alpha)}_\Omega$ when $\alpha\in (0,2)$ were proved quite recently. For the one-dimensional Euclidean case, it was first obtained by Valverde \cite{Valverde2016} by means of probabilistic techniques. Subsequently, using distributional properties of isotropic $\alpha$-stable processes, small time asymptotics of the fractional heat content in higher dimensions was proved by Park and Song \cite{ParkSong2022} for any bounded domain with $C^{1,1}$ regular boundary. More precisely, defining the function $\mu_\alpha$ as 
	\begin{align}\label{eq:mu_alpha}
		\mu_\alpha(t)=\begin{cases}
			t^{1/\alpha}\E[\overline{Y}^{\alpha}(1)] & \mbox{if $1<\alpha\le 2$} \\
			\frac{1}{\pi}t\log(1/t) & \mbox{if $\alpha=1$} \\
			t & \mbox{if $0<\alpha<1$},
		\end{cases}
	\end{align}
	where $\overline{Y}^{\alpha}(t)$ denotes the running supremum of one-dimensional symmetric $\alpha$-stable process,
	it was shown in \cite{ParkSong2022} that 
	\begin{align}\label{eq:alpha_limit_intro}
		\lim_{t\to 0}\frac{|\Omega|-Q^{(\alpha)}_\Omega(t)}{\mu_\alpha(t)}=\begin{cases}
			|\partial\Omega| & \mbox{if $1\le\alpha\le 2$} \\
			\Per_\alpha(\Omega) & \mbox{ if $0<\alpha<1$},
		\end{cases}
	\end{align}
where $|\partial\Omega|$ is the perimeter of $\Omega$, and $\Per_\alpha(\Omega)$ denotes the fractional perimeter of $\Omega$ when $0<\alpha<1$, see \cite[Equation~(1.2)]{ParkSong2022} for details. It is noteworthy that the rate function in the case of $\alpha\in [1,2]$ does not depend on the dimension of the space. This is an indication that $\mu_\alpha$ should be universal if one replaces $\R^n$ by a different space with different geometry.
	
	In this paper, we study the small time asymptotic of the fractional heat content for \emph{Carnot groups}. A connected and simply connected
	 nilpotent Lie group $(\mathbb{G},\star)$ is said to be a \emph{Carnot group of step $k$} if its Lie algebra ${\mathfrak{g}}$  admits a \emph{step $k$ stratification}, that is,
	there exist linear subspaces $V_1,\ldots,V_k$ such that
	\begin{equation}\label{eq:stratification}
		{\mathfrak{g}}=V_1\oplus...\oplus V_k,\quad [V_1,V_i]=V_{i+1},\quad
		V_k\neq\{0\},\quad V_i=\{0\} \ {\textrm{ if }} i>k,
	\end{equation}
	where $[V_1,V_i]$ is the subspace of ${\mathfrak{g}}$ generated by
	the commutators $[X,Y]$ with $X\in V_1$ and $Y\in V_i$. For a good exposition on Carnot groups, we refer to the books \cite{CapognaDanielliPaulsTysonBook,BonfiglioliLanconelliUguzzoniBook}. $(\mathbb G,\star)$ is non-abelian unless $k=1$. Such groups are of special interests as they are equipped with a sub-Riemannian structure, which we are going to introduce now. Set
	$m_i=\dim(V_i)$, for $i=1,\dots,k$ and $n_i=m_1+\dots +m_i$, so that $n_k=N$. For the sake of simplicity, we write
	$n_0=0,\ m:=m_1$.
	We denote by $Q$ the {\em homogeneous dimension} of $\mathbb{G}$, that is, we set
	\begin{align*}
		Q:=\sum_{i=1}^{k} i \dim(V_i).
	\end{align*}
	We choose a basis $\{X_1,\dots,X_N\}$ of
	$\mathfrak{g}$ adapted to the stratification, that is,
	$X_{n_{j-1}+1},\dots,X_{n_j}$ is a basis of $V_j$  for each $j=1,\dots, k$.
	With an abuse of notation, let $X=\{X_1,\dots,X_{N}\}$ be the family
	of left invariant vector fields corresponding to the above basis. The sub-bundle of the tangent bundle $T\mathbb{G}$ that is spanned by the
	vector fields $X_1,\dots,X_{m}$ plays a particularly important
	role in the theory, it is called the {\em horizontal bundle}
	$H\mathbb{G}$; the fibers of $H\mathbb{G}$ are 
	\begin{align*}
	H_x\mathbb{G}=\mbox{span}\{X_1(x),\dots,X_{m}(x)\},\qquad x\in\mathbb{G} .
	\end{align*}
	We can endow each fiber of $H\mathbb{G}$ with an inner 
	product $\langle \cdot,\cdot\rangle$ that makes the basis vectors $X_1(x),\ldots,X_{m}(x)$
	orthonormal. Due to the stratification in \eqref{eq:stratification}, $X_1,\ldots, X_m$ satisfies \emph{H\"ormander's condition}, that is, 
	\begin{align*}
		\mathrm{Lie}(X_1,\ldots, X_m)=T\mathbb{G}.
	\end{align*}
	Since $\mathbb{G}$ is nilpotent, it has a unique (up to multiplicative constants) bi-invariant Haar measure, and $(\mathbb{G},H\mathbb{G})$ is a sub-Riemannian manifold equipped with the smooth bi-invariant Haar measure. Once an orthonormal basis  $X_1,\dots,X_{m}$ of the horizontal bundle is fixed, we
	define the horizontal gradient of a smooth function $f:\mathbb{G}\to \R$, denoted by $\nabla_H f$, as,
	\begin{equation*}
		\nabla_{H}f:=\sum_{i=1}^{m}(X_if)X_i,
	\end{equation*}
	and the associated sub-Laplacian by
	\begin{align*}
		\mathcal{L}f= -\frac{1}{2}\sum_{j=1}^m X_j^2f.
	\end{align*}
	By \cite[Theorem~3.6]{GordinaLaetsch2016} it is known that $\mathcal{L}$ does not depend on the choice of the orthonormal basis of $H\mathbb{G}$. Moreover, by \cite[p.~428]{DriverGrossSaloff-Coste2010}, $(\mathcal{L}, C^\infty_c(\mathbb G))$ is essentially self-adjoint in $L^2(\mathbb G)$, the $L^2$ spaces weighted by the left-invariant Haar measure. Denoting the unique self-adjoint extension of $(\mathcal{L}, C^\infty_c(\mathbb G))$ by $(\mathcal{L},D(\mathcal L))$, we consider the fractional powers defined by means of functional calculus:
	\begin{align*}
		\mathcal{L}_\alpha:=\mathcal{L}^{\frac{\alpha}{2}}=\int_{\sigma(\mathcal{L})} \lambda^{\frac{\alpha}{2}}dE_\lambda,
	\end{align*}
	where $E$ is the spectral measure of $\mathcal{L}$. Given any bounded domain $\Omega$, let us consider the Dirichlet fractional sub-Laplacian $\mathcal{L}^\Omega_\alpha$, which is defined as the Friedrichs extension of $(\mathcal{L}_\alpha, C^\infty_c(\Omega))$. Then,
	 $u_\alpha(t,x)= e^{-t\mathcal{L}^\Omega_\alpha} \mathbbm{1}_\Omega(x)$ solves (in the weak sense) the following parabolic problem with Dirichlet boundary condition:
	\begin{equation}\label{eq:PDE}
	\begin{aligned}
		\left(\frac{\partial}{\partial t}+\mathcal{L}_\alpha\right)u_\alpha(t,x)&=0 \quad \mbox{in $(0,\infty)\times\Omega$} \\
		u_\alpha(t,x)&=0 \quad \mbox{in $(0,\infty)\times\partial\Omega$} \\
		u_\alpha(0,x)&=1 \quad \mbox{for all $x\in\Omega$}.
 	\end{aligned}
 	\end{equation}
 	The fractional heat content of $\Omega$ is defined as
 	\begin{align*}
 		Q^{(\alpha)}_\Omega(t)=\int_\Omega u_\alpha(t,x)dx.
 	\end{align*}
 	
 We also need the following definition of \emph{characteristic points} on the boundary of a domain.
 \begin{definition}[Characteristic points] Let $\Omega$ an open subset of $\mathbb{G}$ with $C^1$ boundary. A point $p\in\partial\Omega$ is said to be a characteristic point if the projection of the outward unit normal $\nu(p)$ onto $H_p\mathbb{G}$ is zero, or equivalently, $H_p\mathbb{G}\subseteq T_p(\partial\Omega)$.
 \end{definition}
 The following is the main result of this article.
\begin{theorem}\label{thm:main}
 		Let $\Omega$ be a bounded, open subset of $\mathbb{G}$ with $C^{2}$ boundary having no characteristic points. Then for all $1\le \alpha\le 2$,
 		\begin{align*}
 			\lim_{t\to 0}\frac{|\Omega|-Q^{(\alpha)}_\Omega(t)}{\mu_\alpha(t)}=|\partial \Omega|_{H},
 		\end{align*}
 		where $\mu_\alpha(t)$ is defined in \eqref{eq:mu_alpha}.
 	\end{theorem}
 	
 	When $\alpha=2$, by a recent work of Rizzi and Rossi \cite{RizziRossi2021}, the asymptotic expansion of $Q^{(2)}_\Omega(t)$ holds on any sub-Riemannian manifold, and $\Omega$ is required to have $C^\infty$ boundary with no characteristic points. In the aforementioned paper, the authors used Savo's technique, see \cite{Savo1999}, adapted to the sub-Riemannian setting, and they obtained asymptotic expansion up to order $5$. Moreover, the first two coefficients in the asymptotic expansion are given in terms of the horizontal perimeter and horizontal mean curvature of the domain at boundary points. We also refer to \cite{CaputoRossi2024}, where the results in \cite{RizziRossi2021} have been generalized to RCD spaces with a local Dirichlet form. However, our framework is different from the aforementioned works as fractional sub-Laplacians are non-local operators, and to the best of our knowledge, heat content asymptotics have not been studied for non-local operators on sub-Riemannian manifolds.
 	
 	The main contribution of our work is to propose a novel approach to study small time asymptotic of fractional heat content on Carnot groups by means of probabilistic and analytic techniques. Using standard subordination argument for fractional powers of generators of Markov processes, see \cite[Section~1.15.9]{BakryGentilLedouxBook}, the fractional heat content can be represented as
 	\begin{align*}
 		Q^{(\alpha)}_\Omega(t)=\int_{\Omega}\P_x\left(\mbox{$B^\alpha(s)\in \Omega$ for all $0\le s\le t$} \right) dx,
 	\end{align*}
 	where $B^\alpha$ is the subordinated horizontal Brownian motion on $\mathbb{G}$, see Section~\ref{sec:HBM} for details. When $\alpha=2$, Tyson and Wang \cite{TysonWangJ2018} used the above probabilistic formula to obtain second order asymptotic expansion of the heat content for $3$--dimensional Heisenberg group. Their method relies on the exit distribution of the horizontal Brownian motion on Heisenberg group along the outward horizontal normal direction of the boundary, as intuitively, the process would most likely exit the domain along this direction. To formalize this heuristic, their proofs involve various technical arguments which are difficult to adapt for subordinated processes on general Carnot groups. We give a simplified argument based on Taylor formula and some refined estimates of the exit probability of horizontal Brownian motion on Carnot group, which leads to the small time asymptotics of the fractional heat content for any $1\le \alpha\le 2$. 
 	
 	Let us now briefly describe the main idea behind our method. We consider a functional version of the fractional heat content, namely, for any nonnegative measurable function $f$ on $\mathbb{G}$, we define 
 	\begin{align}\label{eq:Q_f}
 		Q^{(\alpha)}_f(t)=\int_{\mathbb G}\E_x\left[\inf_{0\le s\le t} f(B^\alpha(s))\right]dx.
 	\end{align}
 	This approach is motivated from a recent work \cite{Sarkar2025_2} of the author on the spectral heat content of isotropic processes in Euclidean space.
 	When $f=\mathbbm{1}_\Omega$, \eqref{eq:Q_f} equals the fractional heat content on $\Omega$. Using the Taylor formula on Carnot groups, thanks to Bonfiglioli \cite{Bonfiglioli2009}, we prove the asymptotics of \eqref{eq:Q_f} for compactly supported smooth functions $f$, and the lower bound of asymptotic heat content follows by an approximation argument. Our method also shows that the lower bound of the asymptotic holds for any bounded open set with finite horizontal perimeter, see Proposition~\ref{prop:liminf}. The regularity conditions on the boundary are needed for proving the upper bound of the limit. We emphasize that our approach mainly requires some estimates of the probability of small exit time of the horizontal Brownian motion (see Theorem~\ref{lem:estimate1}), and a Taylor type expansion of smooth functions on Carnot groups, indicating that this method can be useful for proving heat content asymptotics of non-local operators on more general spaces. In particular, we believe that Theorem~\ref{thm:main} should hold for any sub-Riemannian manifold. At present, we do not know how to prove it.
 	
 	The assumption of non-characteristic boundary is very crucial in our method, and the reason is closely similar to \cite{RizziRossi2021,TysonWangJ2018}. We require $C^2$ regularity of the signed distance function from the boundary of the domain, which ensures that the boundary can be represented as the level set of a $C^2$ function. Even if the domain has $C^\infty$ boundary, in presence of characteristic points, the signed distance function from the boundary may not be Lipschitz continuous, see \cite{AlbanoCannarsaScarinci2018}.
 	
 	The rest of the paper is organized as follows: in Section~\ref{sec:preliminaries} we discuss some properties of the horizontal perimeter and the Taylor formula on Carnot groups, which play significant roles in the proofs. We introduce the horizontal Brownian motion, its subordination, and an estimate of the probability of small exit times in Section~\ref{sec:HBM}. Section~\ref{sec:min_funct} is devoted to proving the asymptotic of the functional defined in \eqref{eq:Q_f}. Finally, Theorem~\ref{thm:main} is proved in Section~\ref{sec:proof_main}.

 	\section{Preliminaries on Carnot group}\label{sec:preliminaries}
 	\subsection{Global coordinates of $\mathbb{G}$}\label{sec:global_coordinates} Since $\mathbb{G}$ is simply connected and nilpotent, the intrinsic exponential map $\Exp:\mathfrak{g}\to\mathbb{G}$ is a global diffeomorphism. As a result, fixing any basis $\{X_1,\ldots, X_N\}$ of the Lie algebra $\mathfrak{g}$, we can identify $\mathbb{G}$ with $\R^N$ using the global coordinates
 	\begin{align*}
 		(x_1,\ldots, x_N) = \Exp\left(\sum_{i=1}^N x_i X_i\right).
 	\end{align*}
 	Due to Baker-Campbell-Dynkin-Hausdorff formula, the group multiplication can be realized on $\R^N$ as
 	\begin{align*}
 		x \star y = \Exp^{-1}\left(\Exp(x) \star \Exp(y)\right), \quad x,y\in\R^N.
 	\end{align*}
 	 With an abuse of notation, we write $\mathbb{G}= (\R^N,\star)$. With this identification, $0\in\R^N$ is the identity element, and the bi-invariant Haar measure on $\mathbb{G}$ coincides with the Lebesgue measure on $\R^N$. Also, there is a natural dilation $\delta_\lambda\in\mathrm{Aut}(\mathbb G)$ defined as
 	\begin{equation}\label{eq:dilation}
 		\delta_\lambda(x_1,\ldots,x_N)=
 		(\lambda \xi_1,\ldots,\lambda^k\xi_k),
 	\end{equation} 
 	where $x=(\xi_1,\dots,\xi_k)\in\R^{m_1}\times\ldots\times \R^{m_k}\equiv\mathbb{G}$.
 	
 	 Indeed, the global coordinates depend on the choice of the basis of $\mathfrak{g}$. A natural way of constructing a basis is as follows: recall that $\dim(V_1)=m$, and fix a basis $\{X_1,\ldots, X_m\}$ of $V_1$ which is orthonormal with respect to the left-invariant inner product on $V_1$. For any multi-index $J=\{j_1,\ldots, j_l\}\subset \{1,2,\ldots, m\}^l$ let us denote the higher order Lie brackets by
 	\begin{align*}
 		X^J=[X_{j_1}[X_{j_2}\cdots]\cdots].
 	\end{align*}
 	Due to the stratification of the Lie algebra $\mathfrak{g}$ in \eqref{eq:stratification}, $X^J\in V_l$ whenever $|J|=l$ and $X^{J_1}, X^{J_2}$ are linearly independent if $|J_1|\neq |J_2|$. Since $\mathrm{span}\{X^J: J\subset \{1,\ldots, m\}\}=\mathfrak{g}$, let $\mathcal{J}\subset \mathcal{P}(\{1,\ldots, m\})$ be such that $\{i\}\in\mathcal{J}$ for $1\le i\le m$, and 
 	\begin{align*}
 		\mathcal{B}=\{X^J: J\in\mathcal{J}\}
 	\end{align*}
 	forms a basis of $\mathfrak{g}$.
 	Throughout the paper, we consider the global coordinates of $\mathbb{G}$ with respect to the basis $\mathcal{B}$ unless stated otherwise.
 	
 	\subsection{Homogeneous norms and distances}We call a function $\|\cdot\|: \mathbb{G}\to [0,\infty)$ a homogeneous norm if
 	\begin{enumerate}
 		\item $\|\delta_\lambda x\|=\lambda \|x\|$ for any $\lambda>0$ and $x\in\mathbb{G}$.
 		\item $\|x\|=0$ if and only if $x=0$.
 	\end{enumerate}
 Any homogeneous norm induces a left-invariant pseudometric on $\mathbb{G}$ defined by
 	\begin{align*}
 		d(x,y)=\|y^{-1}\star x\|, \quad \mbox{$x,y\in\mathbb{G}$}.
 	\end{align*}
 	A pseudometric $d$ defined above is also homogeneous, that is, $d(\delta_\lambda x,\delta_\lambda y)=\lambda d(x,y)$ for all $x,y\in\mathbb{G}$ and $\lambda>0$. Any Carnot group $\mathbb G$ with a sub-Riemannian structure described in the introduction above can be equipped with a Carnot-Carath\'eodory metric defined by
 	\begin{align}
 		d_c(x,y)=\inf\left\{\int_0^1 |\dot{\gamma}(t)|_{H} dt: \gamma(0)=x, \gamma(1)=y, \ \dot{\gamma}(t)\in H_{\gamma(t)}(\mathbb G)\right\}.
 	\end{align}
 	It is known that $(\mathbb{G}, d_c)$ is a path connected metric space, see \cite[Chapter~19]{BonfiglioliLanconelliUguzzoniBook}, and $d_c$ is a homogeneous distance.
 	By \cite[Proposition~5.1.4]{BonfiglioliLanconelliUguzzoniBook}, any two homogeneous norms on $\mathbb{G}$ are equivalent. If $d_\infty$ denotes the homogeneous pseudometric on $\mathbb{G}$ induced by the homogeneous norm
 	\begin{align}\label{eq:d_infinity}
 		\|x\|_\infty=\sup_{1\le j\le k}\{\epsilon_j|(x_{n_{j-1}+1},\ldots, x_{n_{j}})|^{\frac{1}{j}}\}, \quad \epsilon_j>0,
 	\end{align}
 	it was proved in \cite[Theorem~5.1]{FranchiSerapioniCassano2003} that one can choose $\epsilon_1=1$ and $\epsilon_j\in (0,1), j\ge 2$, depending on the group $\mathbb{G}$ so that $d_\infty$ becomes a metric. Therefore, any homogeneous pseudometric $d$ satisfies
 	\begin{align*}
 		c^{-1} d_\infty(x,y)\le d(x,y)\le c d_\infty(x,y) \quad \mbox{for all $x,y\in\mathbb{G}$}
 	\end{align*}
 	for some positive constant $c$ independent of $x,y$.
 	Throughout the paper, we work with an arbitrary homogeneous metric $d$ unless stated otherwise.
 	\subsection{Horizontal perimeter} We start by recalling the definition of functions with bounded variation. For $f\in L^1(\mathbb G)$, the \emph{horizontal variation} of $f$ is defined as
 	\begin{align*}
 		\mathrm{Var}_H(f)=\sup\left\{\sum_{i=1}^m\int_{\mathbb G} fX_i\phi_i dx: \phi_i\in C^\infty_c(\mathbb G), \sum_{i=1}^m |\phi_i|^2\le 1\right\}.
 	\end{align*}
 	$f$ is said to have bounded variation if $\mathrm{Var}_H(f)<\infty$. The space of all functions with bounded variation, denoted by $\BV(\mathbb{G})$, is a Banach space with the norm
 	\begin{align*}
 		\|f\|_{\BV(\mathbb G)}=\|f\|_{L^1(\mathbb G)}+\mathrm{Var}_H(f).
 	\end{align*}
 	Variation is lower semicontinuous, that is, for any $(f_n),f\in L^1(\mathbb G)$ with $\|f_n-f\|_{L^1(\mathbb G)}\to 0$ implies $\liminf_{n\to\infty }\mathrm{Var}_H(f_n)\ge \mathrm{Var}_H(f)$.
 	A measurable set $\Omega\subset \mathbb{G}$ is called \emph{Caccioppoli set} if $\mathbbm{1}_\Omega\in\BV(\mathbb G)$. In this case, the \emph{horizontal perimeter} of $\Omega$ is defined as 
 	\begin{align*}
 		|\partial\Omega|_H:=\mathrm{Var}_H(\mathbbm{1}_\Omega).
 	\end{align*}
 	When $\Omega$ is bounded with $C^1$ boundary, $|\partial\Omega|_H<\infty$, and by \cite[Equation~(3.2)]{CapognaDanielliGarofalo1994}
 	\begin{align}\label{eq:per_hausdorff_1}
 		|\partial\Omega|_H=\int_{\partial\Omega}\left[\sum_{i=1}^m \langle X_i,\nu\rangle^2\right]^{\frac12}d\mathcal{H}^{N-1},
 	\end{align}
 	where $\nu$ is the Euclidean unit outward normal to the boundary, and $\mathcal{H}^{N-1}$ is the $(N-1)$-dimensional Euclidean Hausdorff measure. While \eqref{eq:per_hausdorff_1} provides an explicit relationship between the perimeter and Euclidean Hausdorff measure, it is also possible to represent the horizontal perimeter in terms of the Hausdorff measure with respect to a homogeneous distance on $\mathbb{G}$.
 	Let $d_\infty$ be the homogeneous metric defined in \eqref{eq:d_infinity}, and let $\mathcal{S}^{Q-1}_\infty$ denote the $(Q-1)$-dimensional spherical Hausdorff measure on $\mathbb G$ with respect to that metric, that is,
 	\begin{align}\label{eq:sp_hausdorff}
 		S^{Q-1}_\infty(E)=\lim_{\epsilon\to 0+}\inf\left\{\sum_{i=1}^\infty\frac{\mathrm{diam}(D_{x_i}(t_i))^{Q-1}}{2^{Q-1}}: E\subset \bigcup_{i=1}^\infty D_{x_i}(t_i), \quad t_i\le \epsilon\right\},
 	\end{align}
 	where $D_{x_i}(t_i)=\{x\in\mathbb G: d_\infty(x,x_i)\le t_i\}$ and the diameter is with respect to $d_\infty$.
 	Then by \cite[Proposition~1.9]{Magnani2005} and \cite[Theorem~2.5]{Magnani2005}, we have the following relationship between the horizontal perimeter and $\mathcal{S}^{Q-1}_\infty$.
	\begin{lemma}\label{lem:per_hausdorff}
 		Let $\phi\in C^1(\mathbb G)$, $s\in\R$, and $E=\{x\in\Omega: \phi(x)>s\}$. Assume that $\nabla_H \phi(x)\neq 0$ for all $x\in \partial E$, and $|\partial E|_{H}<\infty$. Then there exists a positive constant $c(\mathbb G)$ depending on $\mathbb G$ such that 
 		\begin{align*}
 			|\partial E|_{H}=c(\mathbb G) S^{Q-1}_\infty(\partial E).
 		\end{align*}
 	\end{lemma}
 	\begin{remark}
 		\cite[Theorem~2.5]{Magnani2005} is stated for $H$-differentiable functions and from the discussion after \cite[Defintion~1.12]{Magnani2005} it follows that any $C^1$ function is $H$-differentiable.
 	\end{remark}
 	Another important fact we need in the proof of Theorem~\ref{thm:main} is the continuity of perimeter of hypersurfaces in $\mathbb{G}$. The following result follows directly from \cite[Theorem~9.1]{DanielliGarofaloNhieu2007}.
 	\begin{lemma}\label{lem:per_cont}
 		Let $\phi\in C^2(\mathbb G)$ be such that $E=\{x: \phi(x)>0\}$ is pre-compact in $\mathbb{G}$ and $|\nabla_{H} \phi|$ is bounded away from $0$ in a neighborhood of $\partial E$. Assume that there exists $r_0>0$ such that for all $r\in [-r_0,r_0]$, the family of domains $E_r=\{x\in\mathbb{G}: \phi(x)>r\}$ have $C^2$ boundary. Then, $\lim_{r\to 0} |\partial E_r|_{H} = |\partial E|_{H}$.
 		
 	\end{lemma}
 	As a consequence of Lemma~\ref{lem:per_hausdorff} and \ref{lem:per_cont}, we get continuity of $(Q-1)$-dimensional $d_\infty$--spherical Hausdorff measure of $C^2$ hypersurfaces, that is, 
 	\begin{align}\label{eq:limit_perimiter}
 		\lim_{r\to 0} \mathcal{S}^{Q-1}_\infty(\partial E_r)=\mathcal{S}^{Q-1}_\infty(\partial E),
 	\end{align}
 	where $E, E_r$ satisfy the conditions of Lemma~\ref{lem:per_cont}.
 	Lastly, we need the coarea formula on $\mathbb{G}$ which is proved in \cite[Corollary~3.6]{Magnani2005}.
 	\begin{lemma}\label{lem:coarea}
 		Let $\phi:\mathbb{G}\to\R$ be a Lipschitz function with respect to the Carnot-Carath\'eodory metric, and $h\in L^1(\mathbb G)$. Then we have,
 		\begin{align*}
 			\int_{\mathbb G} h(x) |\nabla_{H} \phi(x)| dx=c(\mathbb G)\int_{\R}\int_{\phi^{-1}(s)} h(x)d\mathcal{S}^{Q-1}_\infty(x) ds,
 		\end{align*}
 		where $c(\mathbb G)$ is same constant as in Lemma~\ref{lem:per_hausdorff}.
 	\end{lemma}
 	\subsection{Taylor formula} Fix a basis $\{X_1,\ldots, X_N\}$ of the Lie algebra $\mathfrak{g}$. With respect to the global coordinates of $\mathbb{G}$ induced by this basis, we have the following Taylor formula due to Bonfiglioli \cite[Corollary~1]{Bonfiglioli2009}.
 	\begin{theorem} \label{thm:taylor} 
 	For any $\phi\in C^{n+1}(\mathbb{G})$ we have
 		\begin{align*}
 			\phi(x\star h)=\phi(x)&+\sum_{r=1}^n\sum_{i_1,\ldots, i_r=1}^N\frac{h_{i_1}\cdots h_{i_r}}{r!} X_{i_1}\cdots X_{i_r}\phi(x) \\
 			&+\sum_{i_1,\ldots, i_{n+1}=1}^N\frac{h_{i_1}\cdots h_{i_{n+1}}}{n!} \\
 			&\times \int_0^1(X_{i_1}\cdots X_{i_{n+1}})\phi\left(x\star\left(sh_1,sh_2,\ldots, sh_N\right)\right)(1-s)^{n}ds.
 		\end{align*}
 	\end{theorem}
 	We say that $\phi\in C^{1,1}(\mathbb G)$ if $u\in C^1(\mathbb G)$, and there exists $L>0$ such that $|X_i \phi(x)- X_i \phi(y)|\le Ld (x,y)$ for all $1\le i\le N$ and for some homogeneous metric $d$. Applying Theorem~\ref{thm:taylor} we get the following upper bound of the error term in first order Taylor expansion of $C^{1,1}$ functions.
 	\begin{lemma}\label{lem:taylor_bound_0}
 		Suppose that $\phi\in C^{1,1}(\mathbb{G})$ and $d$ is a homogeneous metric on $\mathbb{G}$. Then,
 		\begin{align*}
 			\left|\phi(x\star h)-\phi(x)-\sum_{i=1}^N h_iX_i \phi(x)\right|\le c\sum_{j=1}^kd(h,0)^{j+1}
 		\end{align*}
 		for some constant $c>0$ independent of $x,h$.
 	\end{lemma}
 	\begin{proof}
 		By Theorem~\ref{thm:taylor}, using the first order Taylor expansion, it suffices to show that
 		\begin{align}\label{eq:taylor_0}
 			\sum_{i=1}^N |h_i|\left|\int_0^1 X_i \phi(x\star (sh_1,\ldots, sh_N))ds-X_i\phi(x)\right|\le c \sum_{j=1}^kd(h,0)^{j+1}.
 		\end{align}
 		Using Lipschitz regularity of $X_i \phi$, there exists a constant $L>0$ such that $|X_i \phi(x)-X_i \phi(y)|_{H}\le L d(x,y)$ for all $x,y\in\mathbb G$ and $i=1,\ldots, N$. Therefore, the translation invariance of the distance function and \cite[Lemma~3]{Bonfiglioli2009} yield
 		\begin{align*}
 			\left|\int_0^1 X_i \phi(x\star (sh_1,\ldots, sh_N))ds-X_i\phi(x)\right|\le c_1 d(h,0)
 		\end{align*}
 		for some constant $c_1>0$. On the other hand, by equivalence of homogeneous norms $d$ and $d_\infty$ (see \eqref{eq:d_infinity} for definition) on $\mathbb{G}$ it follows that $|h_i|\le c_2 d(h,0)^j$ if $X_i\in V_j$, $j=1,\ldots, k$ for some constant $c_2$ independent of $h$ and $j$. This proves \eqref{eq:taylor_0} and completes the proof of the lemma.
 	\end{proof}
 	\section{Horizontal Brownian motion and subordination}\label{sec:HBM}
 	With the identification $\mathbb{G}=\R^N$, the horizontal Brownian motion $B=(B(t))_{t\ge 0}$ on $\mathbb{G}$ is the unique solution to the following Stratonovich stochastic differential equation 
 	\begin{align*}
 		dB(t)&=\sum_{i=1}^m X_i(B(t))\circ dW^{(i)}_t \\
 		B(0)&=x,
 	\end{align*}
 	where $W_t=(W^{(1)}_t, \ldots, W^{(m)}_t)$ is the standard Brownian motion in $\R^m$, and $m=\dim(V_1)$. The generator of $B(t)$ on the set of compactly supported smooth functions coincides with the sub-Laplacian $-\mathcal{L}$. In the following, we provide the exact form of $B(t)$ in terms of iterated stochastic integrals of Brownian motion, which was originally proved by Ben Arous \cite{BenArous1989} using the Strichartz formula (see \cite{Strichartz1987}) for deterministic ordinary differential equations on nilpotent Lie groups.
 	
 	For any multi-index $J=\{j_1,\ldots, j_l\}\subset \{1,\ldots, m\}^l$ we write
 	\begin{align*}
 		W^J_t=\int_{T_m(t)} dW^{(j_1)}_{t_1}\cdots dW^{(j_l)}_{t_l}, \quad T_l(t)=\{0\le t_1<\ldots\le t_l\le t\}.
 	\end{align*}
 	For any $l\ge 1$, we denote the permutation group of order $l$ by $\sigma_l$, and for any $\sigma\in\sigma_l$ let us write
 	\begin{align*}
 		e(\sigma)=\mathrm{card}\{1\le j\le l: \sigma(j)>\sigma(j+1)\}.
 	\end{align*}
 	By \cite[Theorem~13]{BenArous1989} there exist real constants $c_J$ such that
 	\begin{align}\label{eq:exponential_rep}
 		B(t)=\Exp(\xi(t))(x), \quad \xi(t)=\sum_{J\in\mathcal{J}}c_JW^J_t X^J
 	\end{align}
 	where $X^J,\mathcal{J}$ are defined in Section~\ref{sec:global_coordinates}, and
 	\begin{align*}
 		W^J_t=\sum_{\sigma\in\sigma_{|J|}}\frac{(-1)^{e(\sigma)}}{|J|^2\dbinom{|J|-1}{e(\sigma)}}W^{J\circ \sigma^{-1}}_t.
 	\end{align*}
 	In particular, $c_i=1$ for all $i=1,\ldots, m$.
 	\begin{notation}\label{not:horizontal}
 		We observe that with respect to the global coordinates of $\mathbb{G}$ induced by the basis $\mathcal{B}$ of $\mathfrak{g}$ we have 
 		\begin{align*}
 			B_i(t)=W^{(i)}_t \quad \text{for all $1\le i\le m$.}
 		\end{align*}
 		We call $B_H(t)=(B_1(t), \ldots, B_m(t))$ the horizontal component of $B(t)$.
 	\end{notation}
 	
 	 The horizontal Brownian motion on $\mathbb G$ is $\frac{1}{2}$-self-similar with respect the dilations $(\delta_\lambda)$ defined in \eqref{eq:dilation}, that is, if $B(0)=0$,
 	\begin{align*}
 		\{\delta_\lambda B_t: t\ge 0\}\overset{d}{=}\{B(t\lambda^2): t\ge 0\} \quad \mbox{for all $\lambda>0$}.
 	\end{align*}
 	
 	For $\alpha\in (0,2]$, consider the $\alpha/2$-stable subordinator $S^\alpha_t$ on $\R$, see \cite[Definition~21.4]{Sato_Book} or \cite[Chapter III.1, p. 73]{BertoinBook}, such that 
 	\begin{align}\label{eq:subordinator}
 		\E\left[e^{-\lambda S^\alpha_t}\right]=e^{-t\lambda^{\frac{\alpha}{2}}} \quad \text{for all $\lambda,t >0$}.
 	\end{align}
 	Assume that $S^\alpha$ is independent of the horizontal Brownian motion $B$.
 	Then, the subordinated process $B^\alpha=(B^\alpha(t))_{t\ge 0}$  defined by
 	\begin{align*}
 		B^\alpha(t):=B(S^\alpha_t)
 	\end{align*}
 	is a left translation invariant Markov process on $\mathbb{G}$, and its generator on the space of compactly supported smooth functions coincides with $-\mathcal{L}_\alpha=-\mathcal{L}^{\alpha/2}$. Due to self-similarity of $S^\alpha$, see Appendix~\ref{sec:appendix}, $B^\alpha$ is also self-similar with index $1/\alpha$, that is, 
 	\[
 	\{\delta_\lambda B^\alpha(t): t\ge 0\}\overset{d}{=} \{B^\alpha(t\lambda^\alpha): t\ge 0\}
 	\] for all $\lambda>0$. When $\alpha=2$, we recover the horizontal Brownian motion. For any bounded domain $\Omega\subset\mathbb{G}$, let us denote the exit time of the subordinated process by
 	\begin{align*}
 		\tau^{(\alpha)}_\Omega=\inf\{t\ge 0: B^\alpha(t)\notin \Omega\}.
 	\end{align*}
 	Then the solution to the fractional Dirichlet problem \eqref{eq:PDE} can be written as
 	\begin{align}\label{eq:u_alpha_P}
 		u_\alpha(x,t)=\P_x(\tau^{(\alpha)}_\Omega>t).
 	\end{align}
 	Due to the above probabilistic formula, the small-time analysis of the heat content $Q^{(\alpha)}_\Omega(t)$ depends on the behavior of the exit probability from the domain as $t\to 0$, which is studied in the following subsection.
 	\subsection{Small time estimates of the exit probability} Let $d$ be any homogeneous metric on $\mathbb{G}$. If $\Omega=D_R(0)$, the $d$--ball of radius $R$ around $0$, then 
 	\[
 	\P_0(\tau^{(2)}_\Omega\le t)=\P_0\left(\sup_{0\le s\le t} d(B(t),0)\ge R\right).\]
 	 The next result provides an upper bound of this exit probability.
 		\begin{theorem}\label{lem:estimate1}
 	There exist constants $c_1, c_2>0$ such that for all $R, t>0$,
 	\begin{equation}\label{eq:exit_probability}
 		\begin{aligned}
 			\P_0\left(\sup_{0\le s\le t}d(B(s),0)> R\right)&\le c_1 e^{-\frac{R^2}{c_2 t}}.
 		\end{aligned}
 	\end{equation}
 	\end{theorem}
 	\begin{remark}
 		We note that by \cite[Proposition~3.3]{Baudoin2025} one has 
 		\begin{align*}
 			\P_0\left(\sup_{0\le s\le t}d(B(s),0)> R\right)\le c_1 \exp\left(-c_2\frac{KR^2}{e^{c_3Kt}-1}\right)
 		\end{align*}
 		for some constants $c_1,c_2, c_3$ and $K\ge 0$ that depends on $\mathbb G$. The above theorem provides sharper upper bound for large values of $t$.
 	\end{remark} 
 	We first prove the following estimate for the tail probability of horizontal Brownian motion.
 	\begin{proposition}\label{prop:estimate2}
 		There exists $c_1,c_2>0$ such that for any $R, t>0$ 
 		\begin{align}\label{eq:first_estimate}
 			\P_0(d(B(t),0)>R)\le c_1 e^{-\frac{R^2}{c_2 t}}.
 		\end{align}
 	\end{proposition}
 	
 	
 	\begin{proof}
 		Since homogeneous metrics on $\mathbb{G}$ are equivalent, we are going to assume that $d=d_c$, the Carnot-Carath\'eodory distance. This would only change the constants in the estimate. Let $p_t$ denote the density of $B(t)$ starting from $0$. By \cite[p.~50]{VaropoulosBook}, there exists a constant $c>0$ such that for all $x\in\mathbb G$ and $t>0$,
 		\begin{align}\label{eq:heat_kernel_est}
 			p_t(x)\le c^{-1} t^{-\frac{Q}{2}} \exp(-d_c(x,0)^2/ct).
 		\end{align}
 		Let $D_R:=D_R(0)$ denote the $d_c$--ball of radius $R$ centered at the origin. Then by \eqref{eq:heat_kernel_est},
 		\begin{align*}
 			\P_0(d_c(B(t),0)>R)\le c^{-1} t^{-\frac{Q}{2}}\int_{D_R^c} e^{-\frac{d_c(x,0)^2}{ct}} dx.
 		\end{align*}
 		Using the co-area formula in Lemma~\ref{lem:coarea} with $\phi(x)=d_c(x,0)$ and noting that $|\nabla_H d_c(x,0)|=1$ for almost every $x$, the above integral reduces to
 		\begin{align}\label{eq:prob_equal}
 			c^{-1} t^{-\frac{Q}{2}}\int_{D_R^c} e^{-\frac{d_c(x,0)^2}{ct}} dx=c(\mathbb G)c^{-1} t^{-\frac Q2} \mathcal{S}_\infty (\partial D_1)\int_R^\infty e^{-\frac{r^2}{ct}} r^{Q-1} dr,
 		\end{align}
 		where $\mathcal{S}_\infty(\partial D_1)$ denote the spherical Hausdorff measure (see \eqref{eq:sp_hausdorff}) of the boundary of the unit ball. In the above equality, we have used the dilation property of $\mathcal{S}_\infty$, that is,
 		\begin{align*}
 			\mathcal{S}_\infty(\delta_\lambda E)=\lambda^{Q-1}\mathcal{S}_\infty(E) \quad \mbox{for all $\lambda>0$ and measurable set $E$}.
 		\end{align*}
 		Absorbing the constants in \eqref{eq:prob_equal} into a single constant and making a change of variable $r/\sqrt{t}\mapsto v$ we get 
 		\begin{align*}
 				\P_0(d_c(B(t),0)>R)\le c'\int_{\frac{R}{\sqrt t}}^\infty e^{-\frac{v^2}{c}} v^{Q-1} dr.
 		\end{align*}
 		As for any $v>0$ there exist universal constants $C,c_2>0$
 		\begin{align*}
 			 e^{-\frac{v^2}{c}} v^{Q-1}\le C e^{-\frac{v^2}{c_2}} v,
 		\end{align*}
 		the above probability can therefore be bounded by
 		\begin{align*}
 				\P_0(d_c(B(t),0)>R)\le c' C \int_{\frac{R}{\sqrt t}}^\infty e^{-\frac{v^2}{c_2}} v=c_1 e^{-\frac{R^2}{c_2 t}}.
 		\end{align*}
 		This completes the proof of the proposition.
 	\end{proof}
 	\begin{proof}[Proof of Theorem~\ref{lem:estimate1}]
 		For notational simplicity let us write 
 		\begin{align*}
 		\tau_R=\inf\{t\ge 0: d(B(t), 0)>R\}.
 		\end{align*}
 		Since $B(t)$ has continuous paths, we have $d(B(\tau_R), 0)=R$. Also, due to the equivalence between homogeneous metrics, let us assume that $d=d_c$, the Carnot-Carath\'eodory distance, as the estimate in \eqref{eq:exit_probability} will change only by a constant. By triangle inequality, we have 
 		\begin{align*}
 			R=d_c(B(\tau_R), 0)\le d_c(B(t),0)+d_c(B(\tau_R), B(t)) \quad \mbox{for all $t\ge 0$}.
 		\end{align*}
 		Therefore we have the following union bound of the exit probability:
 		\begin{align}\label{eq:union_bound1}
 			\P_0(\tau_R\le t)\le \P(d_c(B(t),0)>R/2, \ \tau_R\le t)+\P_0(d_c(B(t), B_{\tau_R})>R/2, \ \tau_R\le t)
 		\end{align}
 		Since $\tau_R$ is a stopping time with respect to the natural filtration generated by $(B(t))_{t\ge 0}$, using the strong Markov property of the horizontal Brownian motion and \cite[Exercise~8.17]{BlumenthalGetoor1968}, the second term above can be bounded as follows:
 		\begin{align*}
 		&\ \  \P_0\left(d_c(B(t), B(\tau_R))>R/2, \ \tau_R\le t\right) \\
 		&=\E_0\left[\mathbbm{1}\{\tau_R\le t\}\P_{B(\tau_R)}(d_c(B(t-\tau_R), B(\tau_R))>R/2)\right] \\
 		&\le \sup_{x\in\mathbb G}\sup_{0\le s\le t} \P_x(d_c(B(s),x)>R/2)) \\
 		& = \sup_{0\le s\le t} \P_0(d_c(B(s),0)>R/2),	
 		\end{align*}
 	where the last identity follows from left-translation invariance of $B(t)$.	By Proposition~\ref{prop:estimate2} it follows that 
 	\begin{align*}
 		\sup_{0\le s\le t} \P_0(d_c(B(s),0)>R/2)\le \sup_{0\le s\le t} c_1e^{-\frac{R^2}{c_2 s}}=c_1e^{-\frac{R^2}{c_2 t}}.
 	\end{align*}
 	Hence, \eqref{eq:union_bound1} yields the desired estimate in \eqref{eq:exit_probability}, which completes the proof of the theorem.
 	\end{proof}

 	\section{Asymptotic of the minimum functional}\label{sec:min_funct}
 	For any nonnegative $f\in L^1(\mathbb G)$, we define the functional of the heat content as
 	\begin{align}
 		Q^{(\alpha)}_f(t)=\int_{\mathbb{G}}\E_x\left[\inf_{0\le s\le t} f(B^\alpha(s))\right]dx.
 	\end{align}
 	Note that the above integral is finite as the right hand side is bounded by $\|f\|_{L^1(\mathbb G)}$. Also, for any $t\ge 0$ we have
 	\begin{align*}
 		Q^{(\alpha)}_{\mathbbm{1}_\Omega}(t)=Q^{(\alpha)}_\Omega(t).
 	\end{align*}
 	In this section, we study the small time asymptotic behavior of $Q^{(\alpha)}_f(t)$ for nonnegative compactly supported smooth functions $f$. We recall that for any $f\in C^\infty_c(\mathbb G)$, we have $\mathrm{Var}_H(f)=\int_{\mathbb G} |\nabla_H f| dx$.

 	\begin{theorem}\label{thm:smooth_function}
 		Let $f\in C^\infty_c(\mathbb G)$ be nonnegative. Then for any $1\le \alpha\le 2$, 
 		\begin{align*}
 			\lim_{t\to 0} \frac{|\Omega|-Q^{(\alpha)}_f(t)}{\mu_\alpha(t)}=\mathrm{Var}_H(f).
 		\end{align*}
 		where $\mu_\alpha$ is defined in \eqref{eq:mu_alpha}.
 	\end{theorem}
 	The proof of this theorem will use Taylor formula on Carnot groups and several small time estimates of the supremum process $\sup_{0\le s\le t} d(B^\alpha(s),0)$, which are proved in subsequent lemmas. We start with a straightforward consequence of Lemma~\ref{lem:taylor_bound_0} applied to functions of left translation invariant processes on $\mathbb{G}$.
 	\begin{lemma}\label{lem:taylor_bound}
 		Fix a homogeneous metric $d$ on $\mathbb{G}$. Let $Z=(Z(t))_{t\ge 0}$ be any left translation invariant process on $\mathbb{G}$ with $Z(0)=0$. For any $x\in\mathbb{G}$ and $f\in C^{1,1}(\mathbb G)$ such that $\|X_i f\|_\infty<\infty$ for all $1\le i\le N$ we have
 		\begin{align*}
 			&  \left|\inf_{0\le s\le t}f(x\star Z(s))-f(x)-\inf_{0\le s\le t}\sum_{i=1}^m Z_i(s) X_i f(x)\right| \\
 			&\le \sum_{i>m}\sup_{0\le s\le t}|Z_i(s)||X_i f(x)|+c\sum_{j=1}^{k-1}\sup_{0\le s\le t} d(Z(s),0)^{1+j} \\
 			&\le c'\sum_{j=1}^{k-1} \sup_{0\le s\le t}d(Z(s),0)^{1+j}
 		\end{align*}
 		for some constant $c'>0$.
 	\end{lemma}
 
\begin{lemma}\label{lem:tail_sup}
	For any $R>0$, as $t\to 0$,
	\begin{align*}
		\P_0\left(\sup_{0\le s\le t}d(B^\alpha(s),0)>R\right)=\mathrm{O}(t).
	\end{align*}
\end{lemma}
\begin{proof}
	Since $S^\alpha=(S^\alpha_t)_{t\ge 0}$ is a process with right continuous paths, for any $t\ge 0$ we have
	\begin{align*}
			\P_0\left(\sup_{0\le s\le t}d(B^\alpha(s),0)>R\right)\le \E_0\left[\P_0\left(\sup_{0\le s\le S^\alpha_t}d(B(s),0)>R\right)\right].
	\end{align*}
Hence, the proof of the lemma follows directly by combining Theorem~\ref{lem:estimate1} and Lemma~\ref{lem:exp_moment}.
\end{proof}
We introduce the following notation that will be used throughout the rest.
\begin{notation}
	Consider the coordinate representation of $B=(B(t))_{t\ge 0}$ according to Notation~\ref{not:horizontal}. We write
	\begin{align*}
		\overline{B}^\alpha_1(t)=\sup\{B^\alpha_1(s): 0\le s\le t\}.
	\end{align*}
	Note that $\overline{B}^\alpha_1(t)$ is the running supremum of symmetric $\alpha$-stable process on $\R$.
\end{notation}
 	
 	 \begin{lemma}\label{lem:sup_tail_limit}
 	  For any $R>0$,
 	 	\begin{align*}
 	 		\lim_{t\to 0}\frac{1}{\mu_\alpha(t)}\E_0\left[\overline{B}^{\alpha}_1(s) \mathbbm{1}\left\{\sup_{0\le s\le t}d(B^{\alpha}(s),0)\le R\right\}\right]=1
 	 	\end{align*}
 	 \end{lemma}
 	 \begin{proof}
 	 	Suppose that $c^{-1} d_\infty \le d\le c d$ for some $c>0$. We note that for any $t,R>0$,
 	 	\begin{align*}
 	 		\left\{\sup_{0\le s\le t}d_\infty(B^{\alpha}(s),0)\le \frac{R}{c}\right\}\subseteq\left\{\sup_{0\le s\le t}d(B^{\alpha}(s),0)\le R\right\}\subseteq \left\{\sup_{0\le s\le t}d_\infty(B^{\alpha}(s),0)\le cR\right\}.
 	 	\end{align*}
 	 	Therefore, it suffices to prove this lemma with $d=d_\infty$. We observe that
 	 	\begin{align*}
 	 	\left\{\sup_{0\le s\le t}d_\infty(B^{\alpha}(s),0)\le R\right\}\subseteq\{\overline{B}^{\alpha}_1(t)\le R\}.
 	 	\end{align*}
 	 	 As a result we get
 	 	\begin{align*}
 	 		0&\le\E_0\left[\overline{B}^{\alpha}_1(t)\mathbbm{1}\left\{\overline{B}^{\alpha}_1(t)\le R\right\}\right]-\E\left[\overline{B}^{\alpha}_1(s) \mathbbm{1}\left\{\sup_{0\le s\le t}d_\infty(B^{\alpha}(s),0)\le R\right\}\right] \\
 	 		&\le R\P_0\left(\sup_{0\le s\le t} d_\infty(B^{\alpha}(s),0)>R\right).
 	 	\end{align*}
 	 	From the exit probability estimate in Lemma~\ref{lem:tail_sup} we have 
 	 	\begin{align*}
 	 		\lim_{t\to 0}\frac{1}{\mu_\alpha(t)}\P_0\left(\sup_{0\le s\le t} d_\infty(B^{\alpha}(s),0)>R\right)=0.
 	 	\end{align*}
 	 	Therefore, it is enough to prove that 
 	 	\begin{align}\label{eq:alpha_limit}
 	 		\lim_{t\to 0}\frac{1}{\mu_\alpha(t)}\E_0\left[\overline{B}^{\alpha}_1(t)\mathbbm{1}\left\{\overline{B}^{\alpha}_1(t)\le R\right\}\right]=1.
 	 	\end{align}
 	 When $1<\alpha\le 2$, $\E[\overline{B}^{\alpha}_1(t)]<\infty$ and by self-similarity of $B^\alpha_1$ we have 
 	 \begin{align*}
 	 	\E\left[\overline{B}^{\alpha}_1(t)\mathbbm{1}\{\overline{B}^{\alpha}_1(t)\le R\}\right]=t^{1/\alpha}\E\left[\overline{B}^{\alpha}_1(1)\mathbbm{1}\{\overline{B}^{\alpha}_1(1)\le R/t\}\right]
 	 \end{align*}
 	 which proves \eqref{eq:alpha_limit}. For $\alpha=1$, we refer to \cite[Proposition~4.3(i)]{Valverde2016}.
 	 \end{proof}
 	 
 	\begin{lemma}\label{lem:tail_expectation}
 		For any $\kappa, R>0$ and $1\le \alpha\le 2$,
 		\begin{align*}
 			\lim_{t\to 0}\frac{1}{\mu_\alpha(t)}\E_0\left[\sup_{0\le s\le t}d(B^\alpha(s),0)^{1+\kappa}\wedge R\right]=0
 		\end{align*}
 		where $\mu_\alpha$ is defined in \eqref{eq:mu_alpha}.
 	\end{lemma}
 	\begin{proof}
 		For any $R>0$ we have 
 		\begin{align}\label{eq:E_kappa}
 			\E_0\left[\sup_{0\le s\le t}d(B^{\alpha}(s),0)^{1+\kappa}\wedge R\right]=(1+\kappa)\int_0^R r^{\kappa}\P_0\left(\sup_{0\le s\le t} d(B^{\alpha}(s),0)>r\right)dr.
 		\end{align}
 	 Also, for any $1\le \alpha\le 2$ we have
 		\begin{align*}
 			\P_0\left(\sup_{0\le s\le t}d(B^{(\alpha)}(s),0)> r\right)&\le \int_0^\infty \P_0\left(\sup_{0\le s\le u} d(B(s),0)>r\right) \P(S^{\alpha}_t\in du).
 		\end{align*}
 		Hence, applying Theorem~\ref{lem:estimate1}, \eqref{eq:E_kappa} yields
 		\begin{equation*}
 		\begin{aligned}
 			\E_0\left[\sup_{0\le s\le t}d(B^{\alpha}(s),0)^{1+\kappa}\wedge R\right]&\le (1+\kappa)c_1\int_0^R r^{\kappa} \left(\E_0\left[\exp\left(-\frac{r^2}{c_2S^\alpha_t}\right)\right]\right)dr
 		\end{aligned}
 		\end{equation*}
 		Using the asymptotic of the integrals in Lemma~\ref{lem:alpha_asymp}, we conclude the proof of this lemma.
 	\end{proof}
 	\begin{proof}[Proof of Theorem~\ref{thm:smooth_function}] Let us denote $S_f=\supp(f)$ and for any $t>0$ we define the following events:
 		\begin{align*}
 			A^x(t)&=\{B^\alpha(s)\in S_f \ \mbox{for all } 0\le s\le t\mid B^\alpha(0)=x \}, \\
 			A^x_r(t)&=\left\{\sup_{0\le s\le t} d_\infty(B^\alpha(s),x)\le r\mid B^\alpha(0)=x\right\},
 		\end{align*}
 		where $d_\infty$ is the metric define in \eqref{eq:d_infinity}.
 		 We note that for any $x\notin S_f$, $f(x)-\E_x\left[\inf_{0\le s\le t} f(B^\alpha(s))\right]=0$, and
 		\begin{align*}
 			\E_x\left[\inf_{0\le s\le t} f(B^\alpha(s))\right]=\E_x\left[\inf_{0\le s\le t}f(B^\alpha(s))\mathbbm{1}_{A^x(t)}\right].
 		\end{align*}
 		
 		Let $R>0$ be large enough such that $S_f\subseteq D_R(0)$, the $d_\infty$--ball of radius $R$ around $0$. Since $B^\alpha$ is left-translation invariant on $\mathbb{G}$, using Lemma~\ref{lem:taylor_bound} for any $f\in C^\infty_c(\mathbb G)$ we have
 		\begin{equation}\label{eq:remainder_bound}
 	\begin{aligned}
 		&\left|\E_x\left[\inf_{0\le s\le t}f(B^\alpha(s))\right]-f(x)-\E_0\left[\inf_{0\le s\le t}\sum_{i=1}^m B^\alpha_i(s) X_i f(x)\mathbbm{1}_{A^0(t)}\right]\right| \\
 		&\le c'\sum_{j=1}^{k-1} \E_0\left[\sup_{0\le s\le t}d_\infty(B^\alpha(s),0)^{1+j}\mathbbm{1}_{A^0(t)}\right] \\
 		&\le  c'\sum_{j=1}^{k-1} \E_0\left[\sup_{0\le s\le t}d_\infty(B^\alpha(s),0)^{1+j}\wedge R^{1+j}\right].
 	\end{aligned}
 	\end{equation}
 	Since $f\in C^\infty_c(\mathbb G)$, invoking Lemma~\ref{lem:tail_expectation} we arrive at 
 	\begin{align}\label{eq:remainder_limit}
 		\lim_{t\to 0}\frac{1}{\mu_\alpha(t)}c\sum_{j=1}^{k-1}\int_{S_f} \E_0\left[\sup_{0\le s\le t}d_\infty(B^\alpha(s),0)^{1+j}\wedge R^{1+j}\right]dx=0.
 	\end{align}
 		We note that
 		\begin{align*}
 			-\E_0\left[\inf_{0\le s\le t}\sum_{i=1}^m B^\alpha_i(s) X_i f(x)\mathbbm{1}_{A^0(t)}\right]=\E_0\left[\sup_{0\le s\le t}\sum_{i=1}^m -B^\alpha_i(s) X_i f(x)\mathbbm{1}_{A^0(t)}\right].
 		\end{align*}
 		For any $x\in S_f$, let us denote 
 		\begin{align*}
 			\delta(x)=\inf\{d_\infty(x,y): y\notin S_f\}.
 		\end{align*}
 		Then, $A^x_{\delta(x)}(t)\subseteq A^x(t)\subseteq A^x_R(t)$ for all $t>0$ and $x\in S_f$. As a result,
 		\begin{equation}\label{eq:two_sided_bounds}
 		\begin{aligned}
 			&\ \ \ \ \E_0\left[\sup_{0\le s\le t}\sum_{i=1}^m -B^\alpha_i(s) X_i f(x)\mathbbm{1}_{A^0_{\delta(x)}(t)}\right] \\
 			&\le \E_0\left[\sup_{0\le s\le t}\sum_{i=1}^m -B^\alpha_i(s) X_i f(x)\mathbbm{1}_{A^0(t)}\right] \\
 			&\le \E_0\left[\sup_{0\le s\le t}\sum_{i=1}^m -B^\alpha_i(s) X_i f(x)\mathbbm{1}_{A^0_{R}(t)}\right].
 		\end{aligned}
 		\end{equation}
 		Since $(B^\alpha_1,\ldots, B^\alpha_m)$ is an isotropic process under the probability law $\P_0$, and for any orthogonal transformation $U:\R^m\to \R^m$, the mapping
 		\begin{align*}
 			x=(\xi_1,\ldots, \xi_k)\mapsto (U\xi_1,\xi_2,\ldots, \xi_k)
 		\end{align*}
 		is a $\|\cdot\|_\infty$--isometry on $\mathbb{G}$, it follows that for any $r>0$,
 		\begin{align*}
 			\E_0\left[\sup_{0\le s\le t}\sum_{i=1}^m -B^\alpha_i(s) X_i f(x)\mathbbm{1}_{A^0_{r}(t)}\right]=|\nabla_{H} f(x)|\E_0\left[\overline{B}^\alpha_1(t)\mathbbm{1}_{A^0_r(t)}\right].
 		\end{align*}
 		Now, combining \eqref{eq:remainder_bound}, \eqref{eq:remainder_limit}, and \eqref{eq:two_sided_bounds} along with Lemma~\ref{lem:sup_tail_limit} yields
 		\begin{align*}
 			& \  \limsup_{t\to 0}\frac{1}{\mu_\alpha(t)}\int_{S_f} \left(f(x)-\E_x\left[\inf_{0\le s\le t} f(B^\alpha(s))\right]\right)dx \\
 			&\le \limsup_{t\to 0}\frac{1}{\mu_\alpha(t)}\E_0\left[\overline{B}^\alpha_1(t)\mathbbm{1}_{A^0_R(t)}\right] \int_{S_f} |\nabla_{H} f(x)|dx \\
 			& =\mathrm{Var}_{H}(f),
 		\end{align*}
 		and
 		\begin{align*}
 			&\  \liminf_{t\to 0}\frac{1}{\mu_\alpha(t)}\int_{S_f} \left(f(x)-\E_x\left[\inf_{0\le s\le t} f(B^\alpha(s))\right]\right)dx \\
 			&\ge \liminf_{t\to 0}\int_{S_f} |\nabla_{H} f(x)| \frac{1}{\mu_\alpha(t)}\E_0\left[\overline{B}^\alpha_1(t)\mathbbm{1}_{A^0_{\delta(x)}(t)}\right] dx \\
 			& \ge \mathrm{Var}_{H}(f),
 		\end{align*}
 		where the last inequality follows from Fatou's lemma. This completes the proof of the theorem.
 	\end{proof}
 	\section{Proof of the lower bound}\label{sec:proof_lower_bound}
 	Let $f$ be a nonnegative measurable function such that $f\in L^1(\mathbb G)$. For any $\epsilon>0$, define the convolution
 	\begin{align*}
 		f_\epsilon(x)=f*\rho_\epsilon(x)=\int_{\mathbb G} f(y^{-1}\star x)\rho_\epsilon(y)dy,
 	\end{align*}
 	where $\rho_\epsilon\in C^\infty_c(\mathbb G)$ is a mollifier supported in $D_\epsilon(0)$, the $d$--ball of radius $\epsilon$ around $0$.
 	\begin{lemma}\label{lem:Q_f}
 		For any $\epsilon>0$, $t>0$ and $0<\alpha\le 2$,
 		\begin{align*}
 			Q^{(\alpha)}_{f_\epsilon}(t)\ge Q^{(\alpha)}_f(t),
 		\end{align*}
 		where $Q_f(t)$ is defined in \eqref{eq:Q_f}.
 	\end{lemma}
 	\begin{proof}
 		By translation invariance of $B^\alpha$ and Fubini's theorem we get 
 		\begin{align*}
 			\E_x\left[\inf_{0\le s\le t} f_\epsilon(B^\alpha(s))\right] &= \E_0\left[\inf_{0\le s\le t} f_\epsilon(x\star B^\alpha(s))\right] \\
 			&=\E_0\left[\inf_{0\le s\le t} \int_{\mathbb G}f(y^{-1}\star x\star B^\alpha(s))\rho_\epsilon(y)dy\right] \\
 			& \ge \E_0\left[\int_{\mathbb G} \inf_{0\le s\le t} f(y^{-1}\star x\star B^\alpha(s))\rho_\epsilon(y) dy\right] \\
 			& = \int_{\mathbb G}\E_{y^{-1}\star x}\left[\inf_{0\le s\le t} f(B^\alpha(s))\right]\rho_\epsilon(y)dy.
 		\end{align*}
 		This shows that for all $\epsilon>0$,
 		\begin{align*}
 			\int_{\mathbb G} \E_{x}\left[\inf_{0\le s\le t} f_\epsilon(B^\alpha(s))\right]dx\ge \int_{\mathbb G} \E_{x}\left[\inf_{0\le s\le t} f(B^\alpha(s))\right]dx.
 		\end{align*}
 		This completes the proof of the lemma.
 	\end{proof}
 	We are now ready to prove the lower bound of the spectral heat content asymptotic.
 	\begin{proposition}\label{prop:liminf}
 		Let $\Omega$ be any bounded, open subset of $\mathbb{G}$ having finite horizontal perimeter. Then,
 		\begin{align*}
 			\liminf_{t\to 0}\frac{|\Omega|-Q^{(\alpha)}_\Omega(t)}{\mu_\alpha(t)}\ge |\partial\Omega|_{H}.
 		\end{align*}
 	\end{proposition}
 	\begin{remark}
 		We note that the lower bound of the limit holds without any regularity assumption on the boundary of the domain.
 	\end{remark}
 	\begin{proof}
 		For any $\epsilon>0$, let us define $f_\epsilon=\mathbbm{1}_\Omega*\rho_\epsilon$, where $\rho_\epsilon\in C^\infty_c(\mathbb G)$ is a mollifier supported in $D_\epsilon(0)$. Also, $\int_{\mathbb G} f_\epsilon(x) dx=|\Omega|$. Therefore by Lemma~\ref{lem:Q_f} and Theorem~\ref{thm:taylor},
 		\begin{align*}
 			\liminf_{t\to 0}\frac{|\Omega|-Q^{(\alpha)}_\Omega(t)}{\mu_\alpha(t)}\ge \liminf_{t\to 0}\frac{\int_{\mathbb G} f_\epsilon(x)dx -Q^{(\alpha)}_{f_\epsilon}(t)}{\mu_\alpha(t)}=\mathrm{Var}_{H}(f_\epsilon).
 		\end{align*}
 	Since $\epsilon>0$ is arbitrary and $f_\epsilon\to \mathbbm{1}_\Omega$ in $L^1(\mathbb G)$, by lower semi-continuity of $\mathrm{Var}_{H}$ we conclude that 
 	\begin{align*}
 			\liminf_{t\to 0}\frac{|\Omega|-Q^{(\alpha)}_\Omega(t)}{\mu_\alpha(t)}\ge\liminf_{\epsilon\to 0} \mathrm{Var}_{H}(f_\epsilon)=|\partial\Omega|_H.
 	\end{align*}
 	This completes the proof of the proposition.
 	\end{proof}
 	\section{Proof of Theorem~\ref{thm:main}}\label{sec:proof_main}
 	Due to Proposition~\ref{prop:liminf}, it suffices to prove that 
 	\begin{align*}
 		\limsup_{t\to 0}\frac{|\Omega|-Q^{(\alpha)}_\Omega(t)}{\mu_\alpha(t)}\le |\partial\Omega|_{H}
 	\end{align*}
 	if $\partial\Omega$ is $C^{2}$ with no characteristic points. In the next lemma, we show that any bounded domain with $C^2$ boundary with no characteristic points can be realized as a hypersurface.
 	\begin{lemma}
 		Let $\Omega$ be a bounded domain with $C^2$ boundary having no characteristic points. Then there exists $\phi\in C^2(\mathbb G)$ such that
 		\begin{enumerate}[(i)]
 			\item $\phi$ is bounded and $|\nabla_{H}\phi|=1$ in a $d_c$--neighborhood of $\partial\Omega$,
 			\item $\Omega=\{x\in\mathbb{G}: \phi(x)>0\}$.
 		\end{enumerate}
 	\end{lemma}
 	\begin{proof}
 		Consider the signed distance function
 		\begin{align*}
 			\delta(x)=\begin{cases}
 				d_c(x,\partial\Omega) & \mbox{if $x\in\Omega$} \\
 				-d_c(x,\partial\Omega) & \mbox{if $x\in\Omega^c$}.
 			\end{cases}
 		\end{align*}
 		Since $\Omega$ has $C^{2}$ boundary with no characteristic points, by \cite[Theorem~1.1]{BossioRizziRossi2025}, there exists $r_0>0$ such that $\delta$ is $C^{2}$ in $\Omega_{r_0}:=\{x\in\mathbb{G}: |\delta(x)|<r_0\}$, and $|\nabla_H \delta|=1$ in $\Omega_{r_0}$. Fix $0<r<r_0$ and define a $C^{2}$ function $\varphi:\R\to\R$ such that 
 		\begin{align*}
 			\varphi(x)=\begin{cases}
 				x & \mbox{if $|x|<r/2$} \\
 				r & \mbox{if $x\ge r$} \\
 				-r & \mbox{if $x\le -r$}.
 			\end{cases}
 		\end{align*}
 	 Consider the function $\phi=\varphi\circ\delta$. As $\delta$ is $C^2$ in $\Omega_{r_0}$ and $\phi$ is constant on $\{x: |\delta(x)|\ge r\}$, $\phi$ is a bounded $C^2$ function on $\mathbb G$. Also, $\phi(x)>0$ if and only if $x\in\Omega$. For any $x\in\Omega_{r/2}$, $|\nabla_H \phi(x)|=|\varphi'(\delta(x)) \nabla_H \delta(x)|=1$. This completes the proof of the lemma.
 	\end{proof}
 	\begin{remark}\label{rem:phi_construction}
 		From the construction of the function $\phi$ in the proof, it follows that $\phi$ satisfies the conditions of Lemma~\ref{lem:per_cont}.
 	\end{remark}
 	We note that for any $t>0$,
 	\begin{align*}
 		|\Omega|-Q^{(\alpha)}_\Omega(t)=\int_{\Omega} \P_x(\tau^{(\alpha)}_\Omega\le t) dx\le \int_{\Omega}\P_x\left(\inf_{0\le s\le t}\phi(B^\alpha(s))\le 0\right)dx.
 	\end{align*}
 	For $\epsilon>0$ let us define $\Omega_\epsilon=\{x\in\Omega: \phi(x)<\epsilon\}$. We will split the integral above into the disjoint sets $\Omega_\epsilon$ and $\Omega\setminus\Omega_\epsilon$, and will prove that the latter term is negligible when scaled by $\mu_\alpha(t)$.
 	\begin{lemma}
 		For any $\epsilon>0$ and $1\le \alpha\le 2$,
 		\begin{align*}
 			\lim_{t\to 0}\frac{1}{\mu_\alpha(t)}\int_{\Omega\setminus\Omega_\epsilon} \P_x\left(\inf_{0\le s\le t}\phi(B^\alpha(s))\le 0\right)dx=0.
 		\end{align*}
 	\end{lemma}
 	\begin{proof}
 		 Assume that $B^\alpha(0)=x$. From Lemma~\ref{lem:taylor_bound} we observe that 
 		\begin{equation}\label{eq:taylor_sup}
 		\begin{aligned}
 			&\left\{\left.\inf_{0\le s\le t}\phi(B^\alpha(s))\le 0\right| B^\alpha(0)=x\right\}\subseteq \left\{\sup_{0\le s\le t}\sum_{i=1}^mB^\alpha_i(s) X_i\phi(x)\right. \\
 			&+\left.\left.c'\sup_{0\le s\le t}\sum_{j=1}^{k-1} d(B^\alpha(s),e)^{1+j}\ge \phi(x)\right| B^\alpha(0)=e\right\}.
 		\end{aligned}
 		\end{equation}
 		From the construction of $\phi$, we have $\|X_i \phi\|_\infty<\infty$ for all $1\le i\le N$. Let 
 		\begin{align}\label{eq:M}
 		M=\max\{\|\nabla_{H}\phi\|_\infty, c'\}.
 		\end{align}
 		 Since $\phi(x)\ge \epsilon$ for any $x\in\Omega\setminus\Omega_\epsilon$, \eqref{eq:taylor_sup} implies that
 		\begin{align*}
 			\sup_{x\in\Omega\setminus\Omega_\epsilon}\P_x\left(\inf_{0\le s\le t}\phi(B^\alpha(s))\le 0\right)&\le \P_0\left(M \overline{B}^\alpha_1(t)\ge \frac{\epsilon}{k}\right)\\
 			&+\sum_{j=1}^{k-1}\P_0\left(M\sup_{0\le s\le t} d(B^\alpha(s),e)^{1+j}\ge \frac{\epsilon}{ck}\right) \\
 			&=: I(t)
 		\end{align*}
 		Since $\overline{B}^\alpha_1(t)\le c \sup_{0\le s\le t}d(B^\alpha(s),0)$ for some constant $c>0$ independent of $t$,
 		by the tail probability estimate in Lemma~\ref{lem:tail_sup}, $I(t)=\mathrm{O}(t)$ as $t\to 0$. Noting that $|\Omega\setminus\Omega_\epsilon|\le |\Omega|$, proof of the lemma follows as $1\le \alpha\le 2$.
 	\end{proof}
 \begin{proof}[Concluding proof of Theorem~\ref{thm:main}]Let $\epsilon>0$ be small enough so that $\phi=\delta$ and $|\nabla_H\phi|=1$ in $\Omega_\epsilon$. Fix $\eta\in (0,1)$. Invoking \eqref{eq:taylor_sup} once again, we have
 	\begin{align*}
 		\int_{\Omega_\epsilon}\P_x\left(\inf_{0\le s\le t}\phi(B^\alpha(s))\le 0\right)dx&\le \int_{\Omega_\epsilon}\P_0\left(|\nabla_H\phi(x)|\overline{B}^\alpha_1(t)\ge \eta\phi(x)\right)dx \\
 		&+\sum_{j=1}^{k-1}\int_{\Omega_\epsilon}\P_0\left(M\sup_{0\le s\le t} d(B^\alpha(s),e)^{1+j}\ge (1-\eta)\frac{\phi(x)}{k-1}\right)dx \\
 		&=:I_1(t)+I_2(t),
 	\end{align*}
 	where $M$ is defined in \eqref{eq:M}. Since $|\nabla_H\phi|=1$ in $\Omega_\epsilon$, by coarea formula in Lemma~\ref{lem:coarea} we obtain
 	\begin{align*}
 		I_1(t)&=c(\mathbb{G})\int_0^\epsilon\int_{\phi^{-1}(r)}\P_x\left(\overline{B}^\alpha_1(t)\ge \eta r\right)d\mathcal{S}^{Q-1}_\infty(x) dr \\
 		&\le c(\mathbb{G})\sup_{0\le r\le \epsilon} \mathcal{S}^{Q-1}_\infty(\phi^{-1}(r))\frac{1}{\eta}\E_0\left[\overline{B}^\alpha_1(t)\wedge \epsilon\right].
 	\end{align*}
 	As $\lim_{t\to 0}\E_0\left[\overline{B}^\alpha_1(t)\wedge \epsilon\right]/\mu_\alpha(t)=1$, we get
 	\begin{align*}
 		\limsup_{t\to 0}\frac{I_1(t)}{\mu_\alpha(t)}\le \frac{c(\mathbb{G})}{\eta}\sup_{0\le r\le \epsilon} \mathcal{S}^{Q-1}_\infty(\phi^{-1}(r)).
 	\end{align*}
 	On the other hand, using coarea formula once again for $I_2(t)$ and following same argument as in $I_1(t)$ we get
 	\begin{align*}
 		I_2(t)\le \frac{1}{1-\eta}\sum_{j=1}^{k-1}\E_0\left[M \sup_{0\le s\le t} d(B^\alpha(s),e)^{1+j}\wedge\frac{\epsilon}{k-1}\right].
 	\end{align*}
 	By Lemma~\ref{lem:tail_expectation}, $\limsup_{t\to 0}I_2(t)/\mu_\alpha(t)=0$ for every $\eta\in (0,1)$ and $\epsilon>0$, which further implies that 
 	\begin{align*}
 		\limsup_{t\to 0}\frac{|\Omega|-Q^{(\alpha)}_\Omega(t)}{\mu_\alpha(t)}\le \frac{c(\mathbb{G})}{\eta}\sup_{0\le r\le \epsilon} \mathcal{S}^{Q-1}_\infty(\phi^{-1}(r)).
 	\end{align*}
 	Letting $\epsilon\to 0$ and $\eta\to 1$ and using \eqref{eq:limit_perimiter}, Lemma~\ref{lem:per_hausdorff} combined with Remark~\ref{rem:phi_construction}, we conclude the proof of the theorem.
 	\end{proof}

\providecommand{\bysame}{\leavevmode\hbox to3em{\hrulefill}\thinspace}
\providecommand{\MR}{\relax\ifhmode\unskip\space\fi MR }
\providecommand{\MRhref}[2]{%
	\href{http://www.ams.org/mathscinet-getitem?mr=#1}{#2}
}
\providecommand{\href}[2]{#2}

	\appendix
	\section{Some estimates related to $\alpha/2$-stable subordinators on $\R$} \label{sec:appendix}
Let $S^\alpha=(S^\alpha_t)_{t\ge 0}$ be the $\alpha/2$-stable subordinator defined in \eqref{eq:subordinator}. Then, $S^\alpha$ is $\frac{\alpha}{2}$-self-similar, that is, 
\begin{align*}
	S^{\alpha}_t \overset{d}{=} t^{\frac{2}{\alpha}}S^{\alpha}_1 \quad \text{for all $t>0$}.
\end{align*}
From \cite[Proposition~28.3]{Sato_Book} it is known that $S^\alpha_t$ is an absolutely continuous random variable for each $t>0$. Denoting the density of $S^\alpha_t$ by $\eta^\alpha_t$, we know from \cite{Skorokhod1961} that
\begin{align*}
	\eta^\alpha_1(u)\sim 2\pi \Gamma(1+\frac{\alpha}{2})\sin(\frac{\pi\alpha}{4}) u^{-1-\frac{\alpha}{2}}  \quad \text{as }u\to\infty,
\end{align*}
and in particular by \cite[p.~97]{BogdanBook}, 
\begin{align}\label{eq:stable_asymptotic}
	\eta^{\alpha}_1(u)\le C\min\{1, u^{-1-\frac{\alpha}{2}}\}.
\end{align}
\begin{lemma}\label{lem:exp_moment}
		For any $0<\alpha<2$ and $\theta_1,\theta_2,c>0$, we have
		\begin{align*}
			\mathbb{E}\left[\exp\left(-\frac{r^{\theta_1}}{c(S^{\alpha}_t)^{\theta_2}}\right)\right]\le C\min\{1, tr^{-\frac{\alpha\theta_1}{2\theta_2}}\}.
		\end{align*}
	\end{lemma}
	\begin{proof}
		From \eqref{eq:stable_asymptotic} we obtain for all $r>0$,
\begin{align*}
			\mathbb{E}\left(\exp\left(-\frac{r^{\theta_1}}{c(S^{\alpha}_t)^{\theta_2}}\right)\right)&=\mathbb{E}\left(\exp\left(-\frac{r^{\theta_1}}{ct^{\frac{2\theta_2}{\alpha}}(S^{\alpha}_1)^{\theta_2}}\right)\right)\\
			&\le C\int_0^\infty \exp\left(-\frac{r^{\theta_1}}{ct^{\frac{2\theta_2}{\alpha}}u^{\theta_2}}\right) u^{-1-\frac{\alpha}{2}}du.
		\end{align*}
		With the change of variable $\frac{r^{\theta_1}}{t^{\frac{2\theta_2}{\alpha}}u^{\theta_2}}\mapsto v$, the above integral simplifies to
		\begin{align*}
			\int_0^\infty \exp\left(-\frac{r^{\theta_1}}{cu^{\theta_2}}\right) u^{-1-\frac{\alpha}{2}}du=\frac{tr^{-\frac{\alpha\theta_1}{2\theta_2}}}{\theta_2}\int_0^\infty e^{-v/c} v^{\frac{\alpha}{2\theta_2}-1}dv.
		\end{align*}
		This completes the proof of the lemma.
	\end{proof}
	
	\begin{lemma}\label{lem:alpha_asymp}
		Let $1\le \alpha\le 2$. Then for $R,c>0$ we have 
		\begin{align}
			\lim_{t\to 0}\mu_\alpha(t)^{-1}\int_0^R r^{\kappa}\E\left[\exp\left(-\frac{r^2}{cS^\alpha_t}\right)\right]dr&=0 \label{eq:A2}
		\end{align}
		where $\mu_\alpha$ is defined in \eqref{eq:mu_alpha}.
	\end{lemma}
	\begin{proof}
		We provide the proof for $1\le \alpha<2$. When $\alpha=2$, $S^\alpha_t=t$ is a deterministic process, and the statement of the lemma follows directly from standard integral estimates. Using Lemma~\ref{lem:exp_moment} with $\theta_1=2, \theta_2=1$, we have
		\begin{align*}
			\int_0^R r^{\kappa}\E\left[\exp\left(-\frac{r^2}{cS^\alpha_t}\right)\right]dr&\le C\int_0^R r^\kappa\min\{1, tr^{-\alpha}\}dr \\
			&= Ct^{\frac{\kappa+1}{\alpha}}\int_0^{Rt^{-\frac1\alpha}}v^{\kappa}\min\{1,v^{-\alpha}\}dv \\
			&=\begin{cases}
				\mathrm{O}(t^{\frac{\kappa+1}{\alpha}}) & \mbox{if $\kappa-\alpha<-1$},\\
				\mathrm{O}(t\log(1/t)) & \mbox{if $\kappa-\alpha=-1$}, \\
				\mathrm{O}(t) & \mbox{if $\kappa-\alpha>-1$}.
			\end{cases}
		\end{align*}
		Since $\kappa>0$, the above upper bound implies \eqref{eq:A2}. This completes the proof of the lemma.
	\end{proof}

	\end{document}